\documentclass[10pt,a4paper]{article}
\usepackage{amsthm}
\usepackage{amscd}
\usepackage{amsmath}
\usepackage{amsthm}
\usepackage{mathrsfs}
\usepackage{latexsym}
\usepackage{amssymb}
\usepackage{amsfonts}
\theoremstyle{plain}

\title{ { \bf{\Large A Poincar\'{e} Lemma for Whitney-de Rham complex }}}

\author{\small  By}

\date{\large  Hou-Yi Chen}

\begin{document}

\maketitle

\begin{center}
 $\mathbf{Abstract}$
 \end{center}
Let $M$ be a real analytic manifold, $Z$ a closed subanalytic subset of $M$.
We show that the Whitney-de Rham complex over $Z$ is quasi-isomorphic to
the constant sheaf $\mathbb{C}_{Z}$.

\tableofcontents

\medskip \noindent
\textbf{}

\section{Introduction}
\label{sec.recall}

In [KS96], Kashiwara-Schapira introduced  the
Whitney functor (real case) and formal cohomology functor (complex case), then
they introduced the notion of ind-sheaves and they also defined Grothendieck six operations in this framework
in   [KS01]. As applications, they defined the Whitney $\mathcal{C}^{\infty}$ functions and Whitney
holomorphic functions on the subanalytic site
as examples of ind-sheaves. The more elementary study for sheaves on the subanalytic site
is performed in [Pr08] and [Pr12].

Let $M$ be a real analytic manifold, by Poincar$\mathrm{\acute{\mathrm{e}}}$ Lemma,
it is well-known that the de Rham complex over $M$ is isomorphic to
$\mathbb{C}_{M}$. The aim of this paper is to show that a theorem of [BP08] follows
easily from a deep result of Kashiwara on regular holonomic $\mathscr{D}$-module [K84] and
the Whitney functor of [KS96]. More precisely,
we show that\\

\noindent
{\bf{Main theorem(=Theorem 3.3.) }} \emph{Let $M$ be a real analytic manifold of dimension $n$
 and $Z$ a closed subanalytic subset of $M$.
Then we have:}
\begin{center}
$\mathbb{C}_{Z}\stackrel{\sim}\longrightarrow
(0\longrightarrow\mathcal{W}_{M,Z}^{\infty}\stackrel{d}\longrightarrow\mathcal{W}_{M,Z}^{(\infty,1)}
\stackrel{d}\longrightarrow\cdots
\stackrel{d}\longrightarrow\mathcal{W}_{M,Z}^{(\infty,n)}\longrightarrow{0})$
\end{center}
\emph{where $\mathcal{W}_{M,Z}^{\infty}$ denotes the sheaf of Whitney functions on $Z$ and
$\mathcal{W}_{M,Z}^{(\infty,i)}$ denotes the sheaf of differential forms
of degree $i$ with coefficients in $\mathcal{W}_{M,Z}^{\infty}$ for each $i$,
i.e., the Whitney-de Rham complex is isomorhpic to $\mathbb{C}_{Z}$.}\\

\noindent
\textbf{Acknowledgments.} I would like to thank Pierre Schapira
for suggesting this problem to me and
for many useful conversations.

\section{Review on Whitney and formal cohomology functors}
\label{sec.recall}

In this section, we review some results on Whitney and
formal cohomology functors. References are
made to [KS96], [KS01], [Pr08], [Pr12] and [S12].

Let $M$ be a real analytic manifold, we denote by $\mathcal{A}_{M}, \mathcal{C}_{M}^{\infty}$
the sheaf of  complex-valued real analytic functions, $\mathcal{C}^{\infty}$-functions on $M$. We denote by
$\mathcal{D}_{M}$ the sheaf of rings on $M$ of finite-order differential operators
with coefficients in $\mathcal{A}_{M}$.

We  denote by $\mathrm{Mod}_{\mathbb{R}-\mathrm{c}}(\mathbb{C}_{M})$ the abelian category of
$\mathbb{R}$-constructible sheaves on $M$ and $\mathrm{Mod}(\mathcal{D}_{M})$
the abelian category of left $ \mathcal{D}_{M}$-modules.
We also denote by $\mathrm{D}^{\mathrm{b}}_{\mathbb{R}-\mathrm{c}}(\mathbb{C}_{M})$ the
bounded derived category consisting of objects whose cohomology groups
belong to   $\mathrm{Mod}_{\mathbb{R}-\mathrm{c}}(\mathbb{C}_{M})$ and
$\mathrm{D}^{\mathrm{b}}(\mathcal{D}_{M})$ the derived category of
$\mathrm{Mod}( \mathcal{D}_{M})$ with bounded cohomologies.\\

\noindent
{\bf{Definition 2.1.}} Let $Z$ be a closed subset of $M$. We denote by $\mathcal{I}_{M,Z}^{\infty}$
the sheaf of $\mathcal{C}^{\infty}$ functions on $M$ vanishing up to infinite
order on $Z$.\\

\noindent
{\bf{Definition 2.2.}} A Whitney function on a closed subset $Z$ of $M$ is an indexed
family $F=(F^{k})_{k\in\mathbb{N}}$ consisting of continuous functions on $Z$
such that  $\forall{m}\in\mathbb{N}$, $\forall{k}\in\mathbb{N}^{n},$ $|{k}|\leq{m}$,$\forall{x}\in{Z}$,
$\forall{\varepsilon}>0$ there exists a neighborhood such that
$\forall{y,z}\in{U\cap{Z}}$
\[ \arrowvert  F^k(z)- \sum_{|j+k|\leq{m}}\frac{(z-y)^{j}}{j!}F^{j+k}(y)\arrowvert
\leq  \varepsilon   d(y,z)^{m-|k|}.  \]
We denote by $W^\infty_{M,Z}$ the space of Whitney
${C}^\infty$ functions on $Z$. We denote by
$\mathcal{W}^\infty_{M,Z}$ the sheaf $U\mapsto W^\infty_{U,U\cap Z}$.

In [KS96], the authors defined the Whitney tensor product functor as follows:
\begin{center}
$\cdot\stackrel{\mathrm{w}}\otimes \mathcal{C}_{M} ^{\infty} \colon
\mathrm{Mod}_{\mathbb{R}-\mathrm{c}}(\mathbb{C}_{M})\rightarrow \mathrm{Mod}( \mathcal{D}_{M})$
\end{center}
in the following way: let $U$ be an open subanalytic subset of $M$
and $Z=M\setminus{U}$. Then $\mathbb{C}_{U}\stackrel{\mathrm{w}}\otimes \mathcal{C}_{M} ^{\infty}=
\mathcal{I}_{M,Z}^{\infty}$ and
$\mathbb{C}_{Z}\stackrel{\mathrm{w}}\otimes \mathcal{C}_{M} ^{\infty}=
\mathcal{W}_{M,Z}^{\infty}$. This functor is exact and extends as a functor in the
derived category, from $\mathrm{D}^{\mathrm{b}}_{\mathbb{R}-\mathrm{c}}(\mathbb{C}_{M})$ to
$\mathrm{D}^{\mathrm{b}}(\mathcal{D}_{M})$.
Moreover the sheaf $F\stackrel{\mathrm{w}}\otimes \mathcal{C}_{M} ^{\infty}$
is soft for any  $\mathbb{R}$-constructible sheaf $F$.

Now let $X$ be a complex manifold and we denote by
$\mathscr{D}_{X}$ the sheaf of rings on $X$ of
finite-order differential operators. We still denote by $X$ the real underlying manifold and
we denote by  $\overline{X}$ the complex manifold conjugate to $X$. One defines the functor of
formal cohomology as follows:

Let $F\in\mathrm{D}_{\mathbb{R}-\mathrm{c}}^{\mathrm{b}}(\mathbb{C}_{X})$, we set

\begin{center}
$F\stackrel{\mathrm{w}}\otimes\mathscr{O}_{X}=R\mathcal{H}om_{\mathscr{D}_{ \overline{X}}}
(\mathscr{O}_{\overline{X}},F\stackrel{\mathrm{w}}\otimes \mathcal{C}_{X}^{\infty})$,
\end{center}
\noindent{where}
$\mathscr{D}_{\overline{X}}$ denotes the
sheaf of rings on
$\overline{X}$ of finite-order differential operators.

Let $M$ be a real analytic manifold, $X$ a complexification of $M$,
$\imath:M\hookrightarrow{X}$ the embedding.
We recall the following result.\\

\noindent
{\bf{Theorem 2.3.}} \emph{([KS96] Theorem 5.10.) Let $F\in\mathrm{D}^{\mathrm{b}}_{\mathbb{R}-c}
(\mathbb{C}_{M})$. Then we have
\begin{center}
$\imath_{*}F\stackrel{\mathrm{w}}\otimes\mathscr{O}_{X}\simeq{\imath}_{*}
(F\stackrel{\mathrm{w}}\otimes\mathcal{C}_{M}^{\infty})$.
\end{center}
In particular,}
\begin{center}
$\mathbb{C}_{M}\stackrel{\mathrm{w}}\otimes \mathcal{O}_{X}\simeq \mathcal{C}_{M}^{\infty}.$
\end{center}
$\hfill \square$

The following proposition is the key point of this note.\\

\noindent
{\bf{Proposition 2.4.}} \emph{([KS96] Corollary 6.2) Let $\mathfrak{M}$ be a regular holonomic $\mathscr{D}_{X}$-module, and let
$F$ be an object of $\mathrm{D}^{\mathrm{b}}_{\mathbb{R}-c}(\mathbb{C}_{X})$.
Then, the natural morphism:
\begin{center}
(2.1)\hfill$R\mathcal{H}om_{\mathscr{D}_{X}}(\mathfrak{M},F\otimes\mathscr{O}_{X})\rightarrow
{R}\mathcal{H}om_{\mathscr{D}_{X}}(\mathfrak{M},F\stackrel{\mathrm{w}}\otimes\mathscr{O}_{X}).\hfill$
\end{center}
is an isomorphism.$\hfill \square$}

\section{Main result}
\label{sec.recall}
Let $X$ be a complex manifold of dimension $n$. We denote by $\mathscr{D}_{X}$ the
sheaf of rings of finite-order differential operators and
$\Theta_{X}$ the sheaf of vector fields on $X$.\\

First we recall the following basic result in $\mathscr{D}_{X}$-module theory.\\

\noindent
{\bf{Proposition 3.1.}} \emph{([K03] Proposition 1.6) The complex}
\begin{center}
$0\rightarrow\mathscr{D}_{X}\otimes_{\mathscr{O}_{X}}\stackrel{n}\bigwedge\Theta_{X}
\rightarrow\cdots\rightarrow\mathscr{D}_{X}\otimes_{\mathscr{O}_{X}}
\stackrel{2}\bigwedge\Theta_{X}\rightarrow\mathscr{D}_{X}\otimes_{\mathscr{O}_{X}}
\Theta_{X}\rightarrow\mathscr{D}_{X}\rightarrow\mathscr{O}_{X}\rightarrow{0}$
\end{center}
\emph{is exact}.$\hfill \square$\\

\noindent
{\bf{Lemma 3.2.}} \emph{Let $\mathscr{M}$ be a left $\mathscr{D}_{X}$-module. Then we have}
\begin{center}
$R\mathcal{H}om_{\mathscr{D}_{X}}(\mathscr{O}_{X},\mathscr{M})\simeq
[\mathscr{M}\rightarrow\Omega_{X}^{1}\otimes_{\mathscr{O}_{X}}
  \mathscr{M}\rightarrow\cdots\rightarrow\stackrel{n}\bigwedge\Omega_{X}^{1}
  \otimes_{\mathscr{O}_{X}}\mathscr{M}].$
\end{center}
\emph{Proof}. By Proposition 3.1, we have
\begin{itemize}
  \item [] $R\mathcal{H}om_{\mathscr{D}_{X}}(\mathscr{O}_{X},\mathscr{M})$
  \item [$\simeq$] $[\mathscr{M}\rightarrow\mathcal{H}om_{\mathscr{D}_{X}}
  (\mathscr{D}_{X}\otimes_{\mathscr{O}_{X}}\Theta_{X},
  \mathscr{M})\rightarrow\cdots\rightarrow\mathcal{H}om_{\mathscr{D}_{X}}
  (\mathscr{D}_{X}\otimes_{\mathscr{O}_{X}}\stackrel{n}\bigwedge\Theta_{X},
  \mathscr{M})]$
  \item [$\simeq$] $[\mathscr{M}\rightarrow\mathcal{H}om_{\mathscr{O}_{X}}
  (\Theta_{X},
  \mathscr{M})\rightarrow\cdots\rightarrow\mathcal{H}om_{\mathscr{O}_{X}}
  (\stackrel{n}\bigwedge\Theta_{X},
  \mathscr{M})]$
  \item [$\simeq$] $[\mathscr{M}\rightarrow\Omega_{X}^{1}\otimes_{\mathscr{O}_{X}}
  \mathscr{M}\rightarrow\cdots\rightarrow\stackrel{n}\bigwedge\Omega_{X}^{1}
  \otimes_{\mathscr{O}_{X}}\mathscr{M}]$
\end{itemize}
where $\Omega_{X}^{1}:=\mathcal{H}om_{\mathscr{O}_{X}}(\Theta_{X},\mathscr{O}_{X})$.
$\hfill \square$\\

Let $M$ be a real analytic manifold, $X$ a complexification of $M$ and $Z$ a closed subanalytic subset of $M$.
We denote by $\Omega_{X}^{1}$ the sheaf of
differential one-form on $X$ and
\begin{itemize}
  \item [(3.1)] $\mathscr{A}_{Z}^{(i)}:=\stackrel{i}\bigwedge\Omega_{X}^{1}\otimes_{\mathscr{O}_{X}}(\mathbb{C}_{Z}\otimes
  \mathscr{O}_{X})$,
  \item [(3.2)] $\mathcal{W}_{M,Z}^{(\infty,i)}:=\stackrel{i}\bigwedge\Omega_{X}^{1}\otimes_{\mathscr{O}_{X}}
  (\mathbb{C}_{Z}\stackrel{\mathrm{w}}\otimes\mathscr{O}_{X})\simeq
  \stackrel{i}\bigwedge\Omega_{X}^{1}\otimes_{\mathscr{O}_{X}}
  (\mathbb{C}_{Z}\stackrel{\mathrm{w}}\otimes\mathcal{C}_{M}^{\infty})
$.
\end{itemize}

Now we are ready to  prove the main theorem of this note below.\\

\noindent
{\bf{Theorem 3.3. }} \emph{Let $M$ be a real analytic manifold of dimension $n$
 and $Z$ a closed subanalytic subset of $M$.
Then we have:}
\begin{center}
$\mathbb{C}_{Z}\stackrel{\sim}\longrightarrow
(0\longrightarrow\mathcal{W}_{M,Z}^{\infty}\stackrel{d}\longrightarrow\mathcal{W}_{M,Z}^{(\infty,1)}
\stackrel{d}\longrightarrow\cdots
\stackrel{d}\longrightarrow\mathcal{W}_{M,Z}^{(\infty,n)}\longrightarrow{0})$
\end{center}
\emph{where $\mathcal{W}_{M,Z}^{\infty}$ denotes the sheaf of Whitney functions on $Z$ and
$\mathcal{W}_{M,Z}^{(\infty,i)}$ denotes the sheaf of differential forms
of degree $i$ with coefficients in $\mathcal{W}_{M,Z}^{\infty}$ for each $i$
which are defined in (3.2),
i.e., the Whitney-de Rham complex is isomorhpic to $\mathbb{C}_{Z}$.}\\

\noindent
{\emph{Proof.}} Take $\mathfrak{M}=\mathscr{O}_{X}$ and $F=\mathbb{C}_{Z}$ in Proposition 2.4.\\
On the one hand, we show that the left hand side of (2.1) is $\mathbb{C}_{Z}$. By Theorem 2.3 and Lemma 3.2,
we get the following complex
\begin{center}
$0\longrightarrow\mathbb{C}_{M}\longrightarrow\mathscr{A}_{M}^{(0)}\stackrel{d}\longrightarrow
\cdots\stackrel{d}\longrightarrow\mathscr{A}_{M}^{(n)}\longrightarrow{0}$
\end{center}
which is exact by Poincar$\acute{\mathrm{e}}$ lemma
where $\mathscr{A}_{M}^{(i)}$'s are defined in (3.1) by taking $Z=M$.
Tensoring $\mathbb{C}_{Z}$, we obtain the following exact sequence
\begin{center}
$0\longrightarrow\mathbb{C}_{Z}\longrightarrow\mathscr{A}_{Z}^{(0)}\stackrel{d}\longrightarrow
\cdots\stackrel{d}\longrightarrow\mathscr{A}_{Z}^{(n)}\longrightarrow{0}$.
\end{center}
Therefore,
\begin{center}
$\mathbb{C}_{Z}\stackrel{\sim}\longrightarrow(0\longrightarrow
\mathscr{A}_{Z}^{(0)}\stackrel{d}\longrightarrow\cdots
\stackrel{d}\longrightarrow\mathscr{A}_{Z}^{(n)}
\longrightarrow{0})$.

\end{center}
On the other hand, the right hand side of (2.1) is the
Whitney-de Rham complex
\begin{center}

$0\longrightarrow\mathcal{W}_{M,Z}^{\infty}\stackrel{d}\longrightarrow\mathcal{W}_{M,Z}^{(\infty,1)}
\stackrel{d}\longrightarrow\cdots
\stackrel{d}\longrightarrow\mathcal{W}_{M,Z}^{(\infty,n)}\longrightarrow{0}$.
\end{center}
Now the result follows from the isomoprhism of (2.1).$\hfill \square$\\

\begin{center}

\end{center}

\noindent
{Institute of Mathematics, Academia Sinica, Taipei 106, Taiwan}\\
E-mail address: houyi@math.sinica.edu.tw

\end{document}